\documentclass[12pt, reqno]{amsart}
\usepackage{amsmath, amstext, amsbsy, amssymb}

\newtheorem{theorem}{Theorem}[section]
\newtheorem{Definition-Lemma}[theorem]{Definition-Lemma}
\newtheorem{lemma}[theorem]{Lemma}
\newtheorem{proposition}[theorem]{Proposition}
\newtheorem{corollary}[theorem]{Corollary}

\theoremstyle{definition}
\newtheorem{definition}[theorem]{Definition}

\theoremstyle{remark}
\newtheorem{remark}[theorem]{Remark}

\numberwithin{equation}{section}

\newcommand{\ad}{\text{ad}}
\newcommand{\be}{\beta}

\newcommand{\C}{ \mathbb C }
\newcommand{\Clg}{ R(\G) }
\newcommand{\CG}{\C [\G]}
\newcommand{\D}{ \mathcal D }
\newcommand{\DA}{{\D}(A)}
\newcommand{\End}{{\rm End}}
\newcommand{\FFock}{{\mathcal F}}
\newcommand{\g}{{\gamma}}
\newcommand{\G}{{\Gamma}}
\newcommand{\Gn}{{\G_n}}
\newcommand{\hD}{\Winf}
\newcommand{\hDA}{\Winf(A)}
\newcommand{\hDg}{\Winf^\G}
\newcommand{\hf}{\frac12}

\newcommand{\Id}{\mbox{Id}}

\newcommand{\la}{\lambda}

\newcommand{\Rz}{R_{\mathbb Z} (\G)}
\newcommand{\pa}{\partial}
\newcommand{\RG}{{\mathcal R}_\G}

\newcommand{\Rgn}{R(\Gn)}
\newcommand{\Tr}{ {\rm Tr} }
\newcommand{\vac}{|0\rangle}

\newcommand{\Winf}{{\mathcal W}_{1+\infty}}
\newcommand{\Z}{ \mathbb Z }

{\vskip-\lastskip\medskip
  \noindent
  {\em #1.}\enspace
  }%
{\qed\par\medskip
  }

\begin{document}

\title[Class algebras of wreath products]
{Vertex algebras and the class algebras of wreath products}

\author[Weiqiang Wang]{Weiqiang Wang}
\address{Department of Mathematics, University of Virginia,
Charlottesville, VA 22904; MSRI, 1000 Centennial Drive,
Berkeley, CA 94720-5070} \email{ww9c@virginia.edu}

\begin{abstract}
The Jucys-Murphy elements for wreath products $\Gn$ associated to
any finite group $\G$ are introduced and they play an important
role in our study on the connections between class algebras of
$\Gn$ for all $n$ and vertex algebras. We construct an action of
(a variant of) the $\Winf$ algebra acting irreducibly on the
direct sum $\RG$ of the class algebras of $\Gn$ for all $n$ in a
group theoretic manner. We establish various relations between
convolution operators using JM elements and Heisenberg algebra
operators acting on $\RG$. As applications, we obtain two distinct
sets of algebra generators for the class algebra of $\Gn$ and
establish various stability results concerning products of
normalized conjugacy classes of $\Gn$ and the power sums of
Jucys-Murphy elements etc. We introduce a stable algebra which
encodes the class algebra structures of $\Gn$ for all $n$, whose
structure constants are shown to be non-negative integers. In the
symmetric group case (i.e. $\G$ is trivial), we recover and
strengthen in a uniform approach various results of
Lascoux-Thibon, Kerov-Olshanski, and Farahat-Higman, etc.
\end{abstract}

\maketitle

\date{}
\section{Introduction}

Let $\G$ be a finite group. Being a generalization of
the symmetric group $S_n$,
the wreath product $\Gn$ is a
semidirect product $\G^n \rtimes S_n$ of the $n$-th product group $\G^n$
and $S_n$ with a natural group structure. In \cite{Wa1},
we realized the direct sum $\RG$ of the representation rings (or
the class algebras) $\Rgn$ of the wreath products $\Gn$ for all $n \ge 0$ as
the Fock space of an infinite-dimensional Heisenberg Lie algebra (cf.
\cite{Mac, Zel} for related algebraic structures on $\RG$; also cf. \cite{FJW}
for further applications). When $\G$ is trivial and thus $\Gn$ becomes
$S_n$, this has been well known (cf. e.g. \cite{Jing}).

As a Fock space carries a canonical structure of
a vertex algebra \cite{Bor, FLM, Kac, FB},
a natural problem is to realize the
vertex algebra structure on $\RG$ in a group-theoretic way. This question was
raised explicitly by I.~Frenkel and the author \cite{FW},
where, as a next step in this direction,
the Virasoro algebra was constructed acting on $\RG$ using certain
convolution products in the class algebras $\Rgn$ for all $n$.
One expects that understanding the connections between $\RG$ and
vertex algebras will unify and shed new lights
on various structures on $\Rgn$ (typically for all $n$).

The goal of the present paper is to develop further
the vertex operator calculus on $\RG$. We construct in a group-theoretic
manner (a variant of)
the $\Winf$ algebra (parametrized by the set $\G^*$ of irreducible characters of $\G$)
acting on $\RG$. The construction makes use of convolution products involving
Jucys-Murphy elements of the wreath products which we introduce in this paper.
As applications of the development of our operator formalism, we obtain
stability properties and algebra generators of the class
algebra $R(\Gn)$. When $\G$ is trivial and thus $\Gn$ becomes $S_n$,
we recover and strengthen in a uniform approach various results
of Lascoux-Thibon, Kerov-Olshanski, and Farahat-Higman, etc,
cf. \cite{LT, KO, FH}.

In  \cite{Wa2} (also cf. \cite{Wa1}), we stressed the significance of parallel (vertex) algebraic structures
behind wreath products, Hilbert schemes, and symmetric products (although
the constructions and proofs can often be different); we
suggested that it is of great value to find the counterparts of
one picture in another
setup, to develop each picture on its own, as well as
to establish connections among different pictures.
(There has been much significant progress
which reveals further similarities and deep connections
among these subjects, cf. e.g. \cite{QW}
for references in these directions.)
The present work can be regarded as another example supporting
this philosophy.
The intuition behind the similarity of the development
here on the class algebras
of wreath products and on the orbifold cohomology rings of symmetric products
in \cite{QW} is
that the study of wreath products can be regarded as the study of the symmetric
products of the orbifold (stack) ${pt}/\G$. While the
definition of Chen-Ruan orbifold
cohomology \cite{CR} is very recent and uses ideas
from the theory of quantum cohomology,
the notion of the class algebra of a finite group is very classical.
As one can put the orbifold $pt/\G$ (resp. $R(\G)$)
on the same footing as a manifold $M$
(resp. cohomology of $M$),
the study of wreath products (rather than the symmetric groups)
remarkably retain many features of the study
of symmetric products. On the other hand, the current work is purely representation-theoretic in
nature and geometry-free.

Let us explain the results and the orgainization of this paper in detail.
In Section~\ref{sec:heis},
we review the basics on wreath products including the description
of conjugacy classes. We recall
our group-theoretic construction of Heisenberg algebra acting on $\RG$
(cf. Theorem~\ref{th:heis}).

In Section~\ref{sec:murphy}, we introduce the Jucys-Murphy (JM) elements
for the wreath product, generalizing the original construction
for the symmetric groups \cite{Juc, Mur}. We establish various
properties of these elements which are analogous to those
of the original JM elements.
We then introduce operators $\mathfrak O(\alpha)$
and $\mathfrak O^k(\alpha)$ in $\End(\RG)$
depending on $\alpha \in \Clg$ linearly and $k \in \Z_+$,
using essentially the convolution product with the $k$-th power
sum of the JM elements for every $n$. The operator $\mathfrak O^1(\alpha)$
goes back to \cite{FW}, while in the symmetric group case
the operator $\mathfrak O^k$ was introduced earlier
in Goulden \cite{Gou} for $k=1$ and in
Lascoux-Thibon \cite{LT} (where the relevance of JM
elements was first pointed out) for general $k$.

We reformulate (cf. Theorem~\ref{th:cubic})
an explicit vertex operator formula in \cite{FW} of $\mathfrak O^1(\alpha)$
in terms of the Heisenberg generators in a way which makes a so-called transfer
property transparent. To do this, we have formulated the class algebra
$\Clg$ as a Frobenius algebra, so it makes sense to talk about the
pushforward map $R(\G) \rightarrow R(\G)^{\otimes n}$ associated to the
diagonal embedding $\G \rightarrow \G^n$, the ``Euler class'' and
the ``Euler number'' of $\G$ (which turns out
to be the number of conjugacy classes in $\G$).
We then  establish a formula (cf. Theorem~\ref{th:covcomm})
involving the operator
$\mathfrak O^k(\alpha)$ and the Heisenberg algebra generators.

In Section~\ref{sec:winf}, we first introduce a family of $\Winf$
algebra parametrized by $\G^*$, which we shall denote by $\hDg$.
We recall that the $\Winf$ algebra is the universal central
extension of differential operators on the circle \cite{Kac}, and
it affords a nice free field realization \cite{Mat, FKRW}. We
realize the action of $\hDg$ on $\RG$ by the commutators of the
operators $\mathfrak O^k(\alpha)$ and the Heisenberg algebra
generators. This result extends the construction of Virasoro
algebra on $\RG$ in \cite{FW} (where only $\mathfrak O^1(\alpha)$
was used) and generalizes the construction in \cite{LT} for the
symmetric groups. To achieve such a realization, we first
establish an explicit identification of the operators $\mathfrak
O^k(\alpha)$ with differential operators on the circle, which
seems new even in the symmetric group case. Our proof uses the
boson-fermion correspondence and the free (fermionic) field
realization of the $\Winf$ algebra. Combining this formula with
Theorem~\ref{th:covcomm}, we derive an explicit formula of
$\mathfrak O(\alpha)$ in terms of a vertex operator. This
specializes to a formula of Lascoux and Thibon \cite{LT} in the
symmetric group case, which was obtained by a very different
method. Such a formula has been used in \cite{LT} to recover
various results scattered in literature. It will be of interest to
apply our general formula for wreath products to derive the
corresponding generalizations.

The operator formalism in this paper is developed in a way
parallel to the development in the orbifold cohomology rings
of symmetric products and cohomology rings of Hilbert schemes of points
on projective surfaces. In particular,
Theorems~\ref{th:heis}, \ref{th:cubic}, and \ref{th:covcomm} have
their exact counterparts in those setups.
The counterparts of the three theorems were
first used by Li, Qin and the author \cite{LQW} to derive
several other results concerning the structure of
the cohomology rings of Hilbert schemes
of points on a surface,
together with some systematic manipulations of
basic constructions such as Heisenberg operators,
pushforwards of diagonals, Euler class of the surface, etc.
This observation allows Qin and the author to axiomatize
the setup which can be applicable for different occasions, including
the study of the orbifold cohomology rings of symmetric products \cite{QW}.
Our current study of the class algebras $\Rgn$ for all $n$ also
fits into the framework, and thus we obtain various further structure
results on $\Rgn$ for free, thanks to
Theorems~\ref{th:heis}, \ref{th:cubic}, and \ref{th:covcomm}
which we have established so far.

Note however that the connections between class algebras of wreath products
and $\mathcal W$ algebras cannot
be derived in the way as in Hilbert schemes and symmetric
products, essentially due to the fact that the square of
the Euler class of a manifold of positive dimension is zero
while any power of the Euler class
of $\G$ is nonzero. Fortunately, we have the tool of
the boson-fermion correspondence available here, which
is nevertheless not effective in the geometric setup.

More explicitly, we obtain stability properties on the convolution
product of normalized conjugacy classes and JM elements of wreath
products, which specialize to the symmetric group case a stability
result stated in Kerov-Olshanski \cite{KO} (cf. \cite{LT} for a
proof). As a consequence, we are able to introduce a {\em stable
algebra} $\mathfrak A_\G$ which encodes the class algebra
structures for $\Gn$ for all $n$. We show that the structure
constants of the algebra $\mathfrak A_\G$ with respect to some
suitable linear basis are nonnegative integers, following the
ideas of Ivanov-Kerov \cite{IK} for the symmetric groups (also cf.
\cite{LT}). Another consequence of our operator formalism is two
distinct sets of algebra generators for the class algebra $\Rgn$.
In the symmetric group case, one set of our algebra generators can
be shown to be equivalent to the one obtained in Farahat-Higman
\cite{FH} in a totally different approach (motivated by its
applications to the block theory). We mention that we have
generalized some other aspects of \cite{FH} to the wreath product
setup and found applications to the cohomology rings of Hilbert
schemes (see \cite{Wa3}).

\noindent {\bf Acknowledgment.}
It is a pleasure to thank Leonard Scott for helpful
discussions on Jucys-Murphy elements. I am grateful
to Jean-Yves Thibon for very helpful email exchange and
for pointing out to me the important reference \cite{IK}
(which readily lead to the proof of a conjecture in
an earlier version on the integrality
and positivity of the structure constants of the
stable algebra $\mathfrak A_\G$).
This work was supported in part by an NSF grant.
Part of the work
is carried out when I enjoy the hospitality of MSRI in
Spring 2002.
\section{The Heisenberg algebra and wreath products}
\label{sec:heis}

\subsection{The class algebra of a finite group}

Given a finite group $\G$, we denote by $\G^*$ the set of complex
irreducible characters and by $\G_*$ the set of conjugacy classes.
Denote by $c^0$ the identity conjugacy class.
We denote by $\Rz$ the $\Z$-span of irreducible characters of $\G$.
We identify $\Clg = \C \otimes_{\Z } R_{\Z}(\G)$ with the space of class
functions on $\G$.

For $c \in \G_*$ let $\zeta_c$ be the order of the centralizer of
an element in the class $c$. Denote by $|\G |$ the order of $\G$.
The usual bilinear form $\langle -, - \rangle_{\G}$ on $R(\G )$ is
defined as follows (often we will omit the subscript $\G$):
\begin{eqnarray}  \label{eq_bilin}
 \langle f, g \rangle = \langle f, g \rangle_{\G}
    = \frac1{ | \G |}\sum_{x \in \Gamma}
          f(x) g(x^{ -1})
 = \sum_{c \in \Gamma_*} \zeta_c^{ -1} f(c) g(c^{ -1}),
\end{eqnarray}
where $c^{ -1}$ denotes the conjugacy class $\{ x^{ -1}, x \in c
\}$. Clearly $\zeta_c = \zeta_{c^{-1}}$. It is well known that
$\langle \g,\g'\rangle =\delta_{\g, \g'},\quad \g, \g' \in \G^*$,
and
\begin{eqnarray}
  \sum_{ \g \in \G^*} \g (c ')  \g ( c^{ -1})
    &= & \delta_{c, c '} \zeta_c, \quad c, c ' \in \G_*. \label{eq_stand}
\end{eqnarray}

One may regard $\CG$ as the space of functions on $\G$, and thus
$\Clg$ as a subspace of $\ad \, \G$-invariants of $\CG$. The {\em
convolution product} in $\CG$ is defined by
\[
  (\be \circ \g) (x) =  \sum_{y \in \G} \be(x y^{-1}) \g(y),
  \quad \be, \g \in \CG, x \in \G.
\]
In particular if $\be, \g \in \Clg$, then so is $\be \circ \g$. We
will often write $\be \circ \g$ as $\be \g$ when there is no
ambiguity. In this way, $\Clg$ is identified with the center of
the group algebra $\CG$, and is often referred to as
the {\em class algebra} of $\G$. It is well known that
\begin{eqnarray}  \label{eq_idem}
   \be\g =\be \circ \g=  \delta_{\be ,\g } h_\g \cdot \g,
    \quad \be, \g \in \G^*,
\end{eqnarray}
where $h_\g =|\G|/d_\g$, and $d_{\g}$ is the {\em degree} of the
irreducible character $\g$. It is known that $h_\g$ is always an
integer.

Given a conjugacy class $c \in \G_*$, we denote by $K^c =\sum_{a
\in c} a$ the characteristic class function of $c$, i.e. the
function which takes value $1$ on elements in the conjugacy class
$c$ and $0$ elsewhere.
\subsection{Conjugacy classes of wreath product $\Gn$}

Given a positive integer $n$, let $\Gamma^n = \Gamma \times \cdots
\times \Gamma$ be the $n$-th direct product of $\Gamma$. The
symmetric group $S_n$ acts on $\Gamma^n$ by permutations: $\sigma
(g_1, \cdots, g_n)
  = (g_{\sigma^{ -1} (1)}, \cdots, g_{\sigma^{ -1} (n)}).
$ The wreath product of $\Gamma$ with $S_n$ is defined to be the
semidirect product
$$
 \Gamma_n = \{(g, \sigma) | g=(g_1, \cdots, g_n)\in {\Gamma}^n,
\sigma\in S_n \}
$$
 with the multiplication
\begin{eqnarray}  \label{eq:twistprod}
 (g, \sigma)\cdot (h, \tau)=(g \, {\sigma} (h), \sigma \tau).
\end{eqnarray}

We denote by $|\la|=\la_1+\cdots+\la_l$ and the length $\ell(\la)
=\ell$, for a partition $\la=(\la_1, \la_2, \cdots, \la_l)$, where
$\la_1 \geq \dots \geq \la_l \geq 1$. We will also make use of
another notation for partitions: $ \la=(1^{m_1}2^{m_2}\cdots) ,$
where $m_i$ is the number of parts in $\la$ equal to $i$. For a
finite set $X$ and $\rho=(\rho(x))_{x \in X}$ a family of
partitions indexed by $X$, we write $\|\rho\|=\sum_{x  \in
X}|\rho(x)|.$ Sometimes it is convenient to regard
$\rho=(\rho(x))_{x \in X}$ as a partition-valued function on $X$.
We denote by ${\mathcal P}(X)$ the set of all partitions indexed
by $X$ and by ${\mathcal P}_n (X)$ the set of all partitions in
${\mathcal P} (X)$ such that $\|\rho\| =n$.

The conjugacy classes of $\Gn$ can be described in the following
way (cf. \cite{Mac}). Let $x=(g, \sigma )\in {\Gamma}_n$, where
$g=(g_1, \ldots , g_n) \in {\Gamma}^n,$ $ \sigma \in S_n$. The
permutation $\sigma $ is written as a product of disjoint cycles.
For each such cycle $y=(i_1 i_2 \cdots i_k)$ the element $g_{i_k}
g_{i_{k -1}} \cdots g_{i_1} \in \Gamma$ is determined up to
conjugacy in $\Gamma$ by $g$ and $y$, and will be called the {\em
cycle-product} of $x$ corresponding to the cycle $y$. For any
conjugacy class $c$ and each integer $i \geq 1$, the number of
$i$-cycles in $\sigma$ whose cycle-product lies in $c$ will be
denoted by $m_i(c)$. Denote by $\rho (c)$ the partition $(1^{m_1
(c)} 2^{m_2 (c)} \ldots )$, $c \in \G_*$. Then each element $x=(g,
\sigma)\in {\Gamma}_n$ gives rise to a partition-valued function
$( \rho (c))_{c \in \G_*} \in {\mathcal P} ( \G_*)$ such that
$\sum_{r, c} r m_r(c) =n$. The partition-valued function $\rho =(
\rho(c))_{ c \in G_*} $ is called the {\em type} of $x$. It is
known that any two elements of ${\Gamma}_n$ are conjugate in
${\Gamma}_n$ if and only if they have the same type. We denote by
$C_\rho$ the conjugacy class of type $\rho$. The corresponding
class function $K^{C_\rho}$ will simply be denoted
by $K^\rho$.

For $\lambda = (1^{m_1} 2^{m_2} \ldots )$, we define
$
  z_{\la } = \prod_{i \geq 1}i^{m_i}m_i!,
$ which is the order of the centralizer of an element of
cycle-type $\la $ in $S_{|\la |}$. The order of the centralizer of
an element $x = (g, \sigma) \in {\Gamma}_n$ of the type $\rho=(
\rho(c))_{ c \in \G_*}$ is
$Z_{\rho}=\prod_{c\in \G_*}z_{\rho(c)}\zeta_c^{l(\rho(c))}.$
%
%
%
%
%
\subsection{The Fock space $\RG$}

Recall that $\Rgn = R_{\Z} (\Gn) \otimes_{\mathbb Z} \mathbb C$.
We set
\begin{eqnarray*}
  \RG = \bigoplus_{n = 0}^\infty R (\Gn).
\end{eqnarray*}
A symmetric bilinear form on $\RG$ is given by
\begin{eqnarray}   \label{eq_series}
\langle u, v \rangle
 = \sum_{ n \geq 0} \langle u_n, v_n \rangle_{\Gn } ,
\end{eqnarray}
where $u = \sum_n u_n$ and $v = \sum_n v_n$ with $u_n, v_n\in
\Gn$.

Denote by $c_n (c \in \G_*)$ the conjugacy class in $\Gn$ of
elements $(x, s) \in \Gn$ such that $s$ is an $n$-cycle and the
cycle product of $x$ lies in the conjugacy class $ c$.
Given $\g \in \Clg$, we denote by $\sigma_n (\g )$ the class
function on $\Gn $ which takes value $n \g (c) $ on elements in
the class $c_n, c \in \G_*$, and $0$ elsewhere.

We define an operator $ \mathfrak p_{ -n} (\gamma)$, where
$ \g \in \Clg, n >0$, to be a
map from $\RG$ to itself by the following composition
\[  R (\G_m) \stackrel{ \sigma_n ( \g ) \otimes }{\longrightarrow}
  \Rgn \bigotimes R (\G_m)  \stackrel{{Ind} }{\longrightarrow}
  R ( {\Gamma}_{n +m}).
\]
We also define an operator $ \mathfrak p_{ n} (\gamma), n
>0$ to be a map from $\RG$ to itself (which is actually the adjoint of
$\mathfrak p_{ -n} (\gamma)$ with respect to the bilinear form
(\ref{eq_series})) as the composition
\[
  R (\G_m)  \stackrel{ Res }{\longrightarrow}
   \Rgn\bigotimes R ( {\G }_{m -n})
 \stackrel{ \langle \sigma_n ( \g), \cdot \rangle }{\longrightarrow}
 R ( {\G }_{m -n}).
\]
We also set $\mathfrak p_0(\g) =0$.
The following theorem was established in \cite{Wa1} in a more
general setup (also see \cite{FJW}).

\begin{theorem} \label{th:heis}
The operators $\mathfrak p_n (\g),$ where $n \in \Z$ and $\g \in \Clg$, satisfy
the Heisenberg algebra commutation relation:
$$
[ \mathfrak p_m (\be), \mathfrak p_n(\g)] = m \delta_{m,-n} \langle \be,\g \rangle.
$$
Furthermore, the space $\RG$ affords an irreducible
highest weight representation of the Heisenberg algebra
with the vacuum vector $\vac = 1 \in \C = R(\G_0)$.
\end{theorem}

In particular, given $\g \in \Clg$ and $y \in R(\G_{n-1})$, we
have by definition that
\begin{eqnarray}   \label{eq:average}
 \mathfrak p_{-1} (\g) (y)
 =\frac1{|\G_{n-1} \times \G|} \sum_{g \in \Gn}\ad\, g (y \otimes\g)
 =\frac1{(n-1)!} \sum_{g \in S_n}\ad\, g (y \otimes\g)
\end{eqnarray}
where we have used $|\G_{n-1} \times \G| =|\G|^{n} (n-1)!$ and
identified $S_n$ as a subgroup of $\Gn$, and the
second equality follows from the fact that $y$ is $\ad\,
\G^{n-1}$-invariant and $\g$ is $\ad \, \G$-invariant.
\section{The Jucys-Murphy elements and convolution operators}
\label{sec:murphy}
\subsection{The Jucys-Murphy elements for wreath products}

Given $c \in \G_*$, we denote by $\kappa_{[21^{n-2}]}^c \in
{\mathcal P}_n (\G_*)$ the function which maps $c$ to the one-part
partition $(2)$, the identity conjugacy class to the partition $( 1^{n-2})$ and other
conjugacy classes to the empty partition. We denote by $K_{[21^{n-2}]}^c$
the conjugacy class in $\Gn$ associated to
$\kappa_{[21^{n-2}]}^c$. By abuse of notation, we also denote by
$K_{[21^{n-2}]}^c \in \Rgn$ the characteristic function
of the conjugacy class $K_{[21^{n-2}]}^c$. In particular, when $c$
is the identity conjugacy class $1$, we will simply write
$\kappa_{[21^{n-2}]}^1$ and $K_{[21^{n-2}]}^1$ as
$\kappa_{[21^{n-2}]}$ and $K_{[21^{n-2}]}$ respectively.
The conjugacy class $K_{[21^{n-2}]}$ can be
decomposed into a disjoint union:

\begin{eqnarray*}
  K_{[21^{n-2}]} = \bigsqcup_{j =1}^n K_{[21^{n-2}]}(j), \qquad K_{[21^{n-2}]}(j) =
  \bigsqcup_{1 \le i<j} K_{[21^{n-2}]}(i,j)
\end{eqnarray*}
where $K_{[21^{n-2}]}(i,j)$ consists of elements $( (g_1, \ldots,
g_n), (i,j)) \in \Gn$, where all $g_k$ except $g_i$ and $g_j$ are
equal to $1\in \G$, and $g_j = g_i^{-1}$  runs over $\G$. It
follows that the cardinality of $K_{[21^{n-2}]}(j)$ is $(j-1)\cdot
|\G| $ and that $K_{[21^{n-2}]}(1)$ is the empty set.

\begin{definition}  \label{def:murphy}
The elements $\xi_j =\xi_{j;n} = \sum_{a \in K_{[21^{n-2}]}(j)} a$
$(j =1, 2, \ldots, n)$ in the group algebra
$\C[\Gn]$ are called the Jucys-Murphy (JM) elements for $\Gn$. Note that $\xi_1=0$.
\end{definition}

When $\G$ is trivial, $(\xi_j)$ are exactly
the JM elements of the symmetric group $S_n$ introduced
in \cite{Juc, Mur}. Our results below indicate that the $(\xi_j)$
in Defintion~\ref{def:murphy} share analogous properties as the
usual Jucys-Murphy elements.

Given a representation $V$ of $\G$ with character $\g \in \Clg$,
the $n$-th outer tensor product $V^{ \otimes  n} $ of $V$
is a representation of the wreath product $\Gn$
whose character will be denoted by $\eta_n ( \g )$: the direct
product $\G^n$ acts on $\g^{\otimes n}$ factor by factor while
$S_n$ by permuting the $n$ factors. Denote by $\varepsilon_n$ the
(1-dimensional) alternating representation of $\Gn$ on which
$\G^n$ acts trivially while $S_n$ acts as the alternating
representation. We denote by $\varepsilon_n ( \g ) \in \Rgn$ the
character of the tensor product of $\varepsilon_n$ and $V^{\otimes
n}$. We extend $\eta_n$ to a map from $\Clg$ to
$\Rgn$ by the following formula (cf. \cite{Wa1}):

\begin{eqnarray*}
  \eta_n (\beta - \g) =
  \sum_{m =0}^n ( -1)^m Ind_{\G_{n -m} \times \G_m }^{\Gn}
   [ \eta_{n -m} (\beta) \otimes \varepsilon_m (\g ) ],
\end{eqnarray*}
where $\beta$ and $\g$ are characters of two $\G$-modules. For
$\g\in \Clg$, the class function $\eta_n(\g) \in \Rgn$ takes value
$\prod_{c \in \G_*} \g(c)^{l(\rho (c))}$ at an element of $\Gn$ of
type $\rho = (\rho(c))_{c \in \G_*}$, cf. \cite{FJW, Mac}.

By construction, $\Gn$ contains $\G^{(i)}$ as a subgroup, where
$\G^{(i)}$ is the $i$-th factor group of $\G^n$. Given $\g \in
\C[\G]$, we denote by $\g^{(i)}$ the corresponding element in $\C[
\G^{(i)}]$, which in turn lies in $\C[\Gn]$.

\begin{theorem}  \label{th:jucysmurphy}
Let $\g \in \Clg$. The $2n$ elements $\xi_j, \g^{(i)},
j=1,2,\ldots,n,$ commute with each other under the convolution
product. Furthermore, We have:

\begin{eqnarray}
\eta_n(\g) &=& (\g^{(1)} +\xi_1)  ( \g^{(2)} +\xi_2)
\cdots   ( \g^{(n)} +\xi_n) ,  \label{eq:trivialchar} \\
\varepsilon_n(\g) &=& (\g^{(1)} -\xi_1)  ( \g^{(2)} -\xi_2) \cdots
( \g^{(n)} -\xi_n) . \nonumber
\end{eqnarray}
\end{theorem}

\begin{proof}

Clearly, $\g^{(k)}$ commutes with $\g^{(j)}$ for $1 \le j,k \le
n$. Note that $\xi_j =\sum_{i<j} \xi_j(i)$ where we have denoted
$\xi_j(i) =\sum_{a \in K_n(i,j)} a$. To check that $\g^{(k)}$
commutes with $\xi_j(i)$ for each $i$ and thus commutes with
$\xi_j$, it suffices to check for $\g = K^c$ with $c\in \G_*$,
which follows from a direct computation.

Next, we shall show that $\xi_j$ and $\xi_\ell$ (say $j <\ell$) commute
with each other. If $k \neq i$ and $k \neq j$, then apparently
$\xi_j(i)$ commutes with $\xi_\ell (k)$. Thus, it suffices to
prove that $\xi_j(i) (\xi_\ell (i) +\xi_\ell(j)) =(\xi_\ell (i)
+\xi_\ell(j)) \xi_j(i).$ This effectively reduces to the case when
$n=3, i=1, j=2, \ell=3$, and $k=1$ or $k=2$, which can be verified
by a simple calculation of multiplications in $\C [\G_3]$:

\begin{eqnarray*}
 \xi_2(1)\, (\xi_3(1) +\xi_3(2))
 &=& \sum_{a,b,c \in \G, abc =1}[ ((a,b,c)  ,(1,3,2)) + ((a,b,c),(1,2,3)) ] \\
 (\xi_3(1) +\xi_3(2)) \, \xi_2(1)
 &=& \sum_{a,b,c \in \G, abc =1} [((a,b,c)  ,(1,2,3)) +
 ((a,b,c),(1,3,2)) ].
\end{eqnarray*}
This finishes the proof that $\xi_j \xi_\ell =\xi_\ell \xi_j$.

It remains to prove the identity (\ref{eq:trivialchar}) since the
proof of the second identity is the same.
We proceed by induction on $n$. When $n=1$, it is clearly true.
Assume that (\ref{eq:trivialchar}) is true for $n-1$, that is,
\begin{eqnarray*}
 (\g^{(1)} +\xi_1)  ( \g^{(2)} +\xi_2)
\cdots    ( \g^{(n-1)} +\xi_{n-1})
 &=& \sum_{\rho \in \mathcal
P_{n-1}(\G_*)}  \prod_{c \in \G_*} \g(c)^{l(\rho (c))} K^{\rho} .
\end{eqnarray*}

It follows from this and the group multiplication
(\ref{eq:twistprod}) in $\Gn$ that
\begin{eqnarray*}
 (\g^{(1)} +\xi_1)  ( \g^{(2)} +\xi_2)
\cdots   ( \g^{(n-1)} +\xi_{n-1}) \g^{(n)}
&=& \sum_{\rho } \prod_{c \in \G_*} \g(c)^{l(\rho (c))} K^{\rho},
\end{eqnarray*}
where $\rho$ runs over the subset of $\mathcal P_{n}(\G_*) $ which
consists of elements obtained from $\mathcal P_{n-1}(\G_*) $ by
all possible ways of adding a disjoint one-cycle.

On the other hand, we observe that, for $\rho \in \mathcal
P_{n-1}(\G_*)$, the product $K^{\rho} \circ \xi_n$ is the sum of
characteristic class function $K^{\rho'}$, where $\rho'$ runs over
elements in $\mathcal P_{n}(\G_*)$ obtained from $\rho$ by all
possible ways of changing a part equal to $i$ to a part equal to
$(i+1)$.

Combining the above observations together, we have established
(\ref{eq:trivialchar}).
\end{proof}

\begin{corollary}
Given $c \in \G_*$, we have
\begin{eqnarray*}
 \prod_{j =1}^n
 \left(\xi_j+ K^{c(j)} \right) = \sum_{\rho \in \mathcal P_n(c)}
 K^\rho,
\end{eqnarray*}
where $\mathcal P_n(c)$ denotes the set of all the
partition-valued functions $\rho \in \mathcal P_n(\G_*)$
such that $\rho (c') =\emptyset$ for any $c' \neq c$. Also, we have

\begin{eqnarray}   \label{eq:trivial}
 \prod_{j =1}^n
 \left(\xi_j+ \sum_{g \in \G^{(j)}} g \right) = \sum_{a \in \Gn} a .
\end{eqnarray}
\end{corollary}

\begin{proof}
The first part is a special case of Theorem~\ref{th:jucysmurphy}
for $\g =K^c$.

Set $\g =\sum_{g\in \G} g$ to be the trivial character of $\G$ in
Theorem~\ref{th:jucysmurphy}. Observe that $\eta_n(\g)$ is the
trivial character of $\Gn$, and thus is equal to $\sum_{a \in \Gn}
a$. Therefore the identity (\ref{eq:trivial}) follows from
Theorem~\ref{th:jucysmurphy}.

\end{proof}

\begin{remark}
A set of coset representatives for $\Gn/\G_{n-1}$ consists of the
elements $g \in \G^{(n)}$ and $a \in K_{[21^{n-2}]}(n)$ which
appear in the last factor of the product in (\ref{eq:trivial}).
The identity (\ref{eq:trivial}) follows also from this fact by
induction.
\end{remark}
\subsection{Frobenius algebra structure on $\Clg$}

 The characteristic function of the identity
in the group $\G$ is the identity of the class algebra $\Clg$, and clearly,
 $$1 =\sum_{\g \in \G^*}  \g/ h_\g. $$
We define a trace map $\Tr: \Clg \rightarrow \C$ by letting $\Tr
(K^c) = \delta_{c,1}/|\G|, c \in \G_*$, and then extending
linearly to $\Clg$. It follows that $\Tr (\g) = h_\g^{-1}$
for $\g \in\G^*$. This induces a bilinear form on $\Clg$ by
$\langle \alpha,\be \rangle =\Tr(\alpha \circ \be)$. One has that $\langle
\alpha,\be \rangle =\delta_{\alpha,\be}$ for $ \alpha,\be \in \G^*$, so this
bilinear form coincides with the standard
one on $\Clg$. In other word,
we have endowed a {\em Frobenius algebra} structure on $\Clg$.

We define a pushforward map $\tau_{2*}: \Clg \rightarrow \Clg
\otimes \Clg$ to be the linear map adjoint to the convolution
product with respect to the bilinear form $\langle -,-\rangle$. We
further define $\tau_{k*}: \Clg \rightarrow \Clg^{\otimes k}$
inductively by the composition

$$\tau_{(k+1)*} = (\tau_{2*} \otimes \mbox{Id}^{\otimes (k-1)})\cdot
\tau_{k*}.$$

\begin{lemma} \label{lem:push}
 The pushforward map $\tau_{2*}: \Clg \rightarrow
\Clg \otimes \Clg$ is characterized by:
 $$\tau_{2*} (\g) =h_\g \g \otimes \g, \quad \g \in\G^*.$$
\end{lemma}

\begin{proof}
Since $\tau_{2*}$ is a linear map, it is characterized by its
values over a linear basis of $\Clg$. For $\alpha, \be, \g \in
\G^*$, we have by definition that
 $$\langle \tau_{2*}(\g), \alpha \otimes \be \rangle_{\G \times \G}
 =\langle \g, \alpha \circ \beta\rangle
 =\langle \g, h_\be  \delta_{\alpha,\be} \be\rangle
 = {h_\be} \delta_{\alpha,\be} \delta_{\g, \be} $$

 On the other hand, we have
 $$\langle {h_\g}  \g \otimes \g, \alpha \otimes \beta \rangle_{\G \times \G}
  = {h_\g}  \delta_{\g, \alpha} \delta_{\g, \be}.$$
 The lemma follows from the comparison of these two identities.
\end{proof}

\begin{remark}
Since $\tau_{2*}(1) =\sum_{\g \in\G^*} \g \otimes \g$,
we find that the composition map of convolution
product and the pushforward $\tau_{2*}$
$$
 \Clg \stackrel{\tau_{2*}}{\longrightarrow} \Clg \times \Clg
 \stackrel{\circ}{\longrightarrow} \Clg $$
sends $1 \in \Clg$ to $\chi:= \sum_{\g\in \G^*} h_{\g} \g$, whose trace is
given by $\Tr(\sum_{\g\in \G^*} h_{\g} \g) =|\G^*|$.
In analog with geometry, the class function
$\chi= \sum_{\g\in \G^*} h_{\g} \g  $ and the number $|\G^*|$
should be regarded as the ``Euler class" and respectively the
``Euler number" of $\G$.
\end{remark}
\subsection{Definition of convolution operators}

\begin{lemma}
For a given $\alpha \in \Clg$, any symmetric function in the $n$ elements
$\xi_1+\alpha^{(1)}, \ldots, \xi_n+\alpha^{(n)}$ lies in $\Rgn$.
\end{lemma}

\begin{proof}
We denoted by $E_i(\g)$ the $i$-th
elementary symmetric functions ($i=1, \ldots,n$) in $\g^{(1)} +\xi_1, \ldots,
\g^{(n)} +\xi_n$.
For $c\in \G_*$, $E_i(K^c)$ is the sum of all elements in $\Gn$ which satisfy the following
property: it has $n-k$ cycles, $(i-k)$ of which is of cycle type $c$ while
$(n-i)$ of which is of cycle type $1$, where $k$ runs over $0,1,\ldots, i$.
This can be seen by directly multiplying them out term by term.
(By the way, this generalizes the following well-known fact about JM elements in
$S_n$: the $i$-th elementary symmetric function of JM elements
is the sum of all permutations having exactly $n-i$ cycles.)
In particular, $E_i(K^c)$ is a class function of $\Gn$. We can see similarly that
$E_i(\g)$, where $\g \in \Clg$, is a class function of $\Gn$.
Any symmetric function in $\g^{(1)} +\xi_1, \ldots,
\g^{(n)} +\xi_n$ can be written as a polynomial in the $E_i(\g)$'s,
and thus is a class function of $\Gn$.
\end{proof}

\begin{corollary}
For $\alpha \in \Clg$, the sum $\sum_{i=1}^n  \xi_i^{k} \circ \alpha^{(i)}$ lies in $\Rgn$.
\end{corollary}

\begin{proof}
The above lemma implies that
 $\sum_{i=1}^n  e^{(\g^{(i)} + \xi_i)} = \sum_{i=1}^n e^{\xi_i} \circ (e^{\g})^{(i)}$
is a class function for each $\g \in \Clg$. As $e^{\g}$ span $\Clg$
when $\g$ ranges over $\Clg$, the corollary follows.
\end{proof}

\begin{definition}
We define $\Xi_n^k$ to be the $k$-th power
sum $\xi_1^k +\ldots +\xi_n^k$ of the JM elements in $\Gn$.
For  $\alpha \in \Clg$, we define the class function $\Xi_n^k(\alpha) \in
\Rgn$ to be
 $$\Xi_n^k(\alpha) =\sum_{i=1}^n  \xi_i^{k} \circ \alpha^{(i)}.$$
We define $\Xi_n(\alpha) =\sum_{i=1}^n e^{\xi_i}\circ
\alpha^{(i)},$ and $\Xi_{n,\hbar}(\alpha) = \sum_{k\ge 0}
{\hbar^k}/{k!} \cdot \Xi_n^k(\alpha)$, where
$\hbar$ is a formal parameter. We further define the
operator $\mathfrak O^k(\alpha) \in \End (\RG)$ (resp. $\mathfrak
O (\alpha)$, or $\mathfrak O_{\hbar}(\alpha)$) to be the
convolution product with $\Xi_n^k(\alpha)$ (resp. $\Xi_n(\alpha)$,
or $\Xi_{n,\hbar}(\alpha)$) in $\Rgn$ for every $n \ge 0$.
\end{definition}

We observe that
the class function $\Xi_n^1(K^c)$ coincides with the characteristic
function $K_{[21^{n-2}]}^c$ of the conjugacy class in $\Gn$ of
type $\kappa_{[21^{n-2}]}^c$.
The operator $\mathfrak O^1(K^c)$ was introduced in \cite{FW} and it
will continue to play a distinguished role in this paper. We
shall denote $\mathfrak O^1(\alpha)$ by $\mathfrak b(\alpha)$ for any $\alpha \in \Clg$, and set
$\mathfrak b = \mathfrak b (1)$.
\subsection{A new formulation of results in \cite{FW}}

We define $:\mathfrak p_m(\alpha_1) \mathfrak p_n(\alpha_2):$ to
be $\mathfrak p_m(\alpha_1) \mathfrak p_n(\alpha_2) $ if $m \le n$
and $ \mathfrak p_n(\alpha_2) \mathfrak p_m(\alpha_1)$ if $m >n.$
We define $$:\mathfrak p(\alpha_1) \mathfrak p (\alpha_2)(z): =
\sum_{m,n\in \Z} :\mathfrak p_m(\alpha_1) \mathfrak p_n(\alpha_2):
z^{-m-n-2}, $$
and similarly define inductively $:\mathfrak p(\alpha_1) \ldots
\mathfrak p (\alpha_k)(z):$ $(k\ge 2)$ from the right to the left
as usually done in the theory of vertex algebras (cf. e.g. \cite{Kac}).

If we write $\tau_{k*}(\alpha) =\sum_i \alpha_{i,1} \otimes \ldots
\otimes \alpha_{i,k}$, then we define the following vertex
operator

$$:\mathfrak p(z)^k:(\tau_{k*} \alpha)
=\sum_i:\mathfrak p(\alpha_{i,1})(z) \ldots \mathfrak p
(\alpha_{i,k})(z):,
$$
and its components by $:\mathfrak p(z)^k:(\tau_{k*} \alpha)
=\sum_{n \in\Z}:\mathfrak p^k:_n(\tau_{k*} \alpha)z^{-n-k}.$ In
other words, the operator $:\mathfrak p^k:_n(\tau_{k*} \alpha)$
lies in $\End(\RG)$.

\begin{definition} \label{def:vir}
 We define $L( \be)(z) = \sum_{n\in \Z} L_n(\be) z^{-n-2} =\hf
:\mathfrak p(z)^2: (\tau_{2*}\be)$ for $\be \in \Clg$.

We say that the commutator of two families of operators $A(\alpha)$ and $
B(\be)$ (where $\alpha,\be\in \Clg$) satisfies the {\em transfer
property} if $[A(\alpha), B(\be)] =[A(\alpha\be), B(1)] =[A(1),
B(\alpha\be)].$
\end{definition}

\begin{theorem}  \label{th:virasoro}
The commutator between $\mathfrak b(\be)$ and the Heisenberg
algebra generator $\mathfrak{p}_{n} (\g)$ is given by
 \[
 [ \mathfrak b(\be), \mathfrak{p}_{n} (\g)] = -n L_n(\be \g),\quad \be,\g \in \Clg,
\]
where the operators $L_n$ acting on the space $\RG$ satisfy the
Virasoro commutation relations:
 \[
    [ L_n(\be), L_m(\g) ]
    = (n-m) L_{n +m}(\be \g) +
     \frac{n^3 -n}{12} \delta_{n, -m} {\Tr}(\chi \be \g) \cdot \text{Id}_{\RG}.
 \]
In particular, the transfer property holds for both commutators.
\end{theorem}

\begin{proof}
To verify the Virasoro relation, it suffices to do so for
$\be,\g\in \G^*$, since all the terms involved depends on $\be$ and
$\g$ linearly. Since $[\mathfrak p_m(\be),\mathfrak p_n(\g)] =m
\delta_{m,-n}\delta_{\be,\g}$, the components $L_n^{\be}$ of
$L^\be (z)\stackrel{\rm def}{=}\hf :\mathfrak p(\be)(z)^2:$ for $
\be \in \G^*$ (as operators acting on the space $\RG$) satisfy the
commutation relation:
 \[
    [L_n^{\be}, L_m^{\g}]
    = \delta_{\be,\g} (n-m) L_{n +m}^{\g} +
     \frac{n^3 -n}{12} \delta_{\be,\g} \delta_{n, -m}.
 \]

Since $\tau_{2*}\g =h_\g  \g \otimes \g$, we have by
Definition~\ref{def:vir} and Lemma~\ref{lem:push} that $L_n(\g)
=h_\g L_n^\g$. We also calculate
$$\Tr(\chi \be \g) = \Tr(\delta_{\be,\g} h_\g^3 \g) = \delta_{\be,\g} h_\g^2.$$
Thus, we obtain that

\begin{eqnarray*}
 [L_n(\be), L_m(\g)] &=&\delta_{\be,\g} h_\be^2 \;[L_n^{\be},L_m^{\g}] \\
 &=& \delta_{\be,\g} h_\be^2 (n-m) \cdot L_{n +m}^{\be} +
  \frac{n^3 -n}{12}\; h_\be^2 \delta_{\be,\g} \delta_{n, -m}  \cdot \text{Id}_{\RG}\\
 &=&(n-m) L_{n +m}(\be  \g) +
     \frac{n^3 -n}{12} \delta_{n, -m} {\Tr}(\chi \be \g) \cdot \text{Id}_{\RG}.
\end{eqnarray*}

To prove the first commutation relation in the theorem, it
suffices to do so for $\g\in \G^*$ and $\be =K^c$ for $c\in \G_*$,
since both sides of the commutation relation depends linearly on
$\be$ and $\g$. This can be checked by using Theorem~4 in
\cite{FW} as follows. Note that $\Delta_c$ and $ \tilde{a}_n(\g)$
therein are $\mathfrak b (K^c)$ and $ \mathfrak p_n(\g)$ respectively in
our current notation.
 \begin{eqnarray*}
 [\mathfrak b (K^c), \mathfrak p_n(\g)]
 &=& - {n \zeta_c^{-1} h_\g^2 \g(c^{-1})}\, L_n^\g \\
 &=& - n \zeta_c^{-1} h_\g \g(c^{-1}) \, L_n (\g)
 = -n \,L_n(K^c \g),
 \end{eqnarray*}
where we have used
 $$K^c \g = \sum_{\be \in \G^*}
\zeta_c^{-1}
 \be(c^{-1})\; \be \g = \zeta_c^{-1} h_\g \g(c^{-1})\, \g.$$
\end{proof}

\begin{remark}
The operators $L_n(1)$ generate the standard Virasoro algebra:
 \[
    [ L_n(1), L_m(1) ]
    = (n-m) L_{n +m}(1) +
     \frac{n^3 -n}{12} \delta_{n, -m} |\G^*| \cdot \text{Id}_{\RG}.
 \]
Note that the central charge for the Virasoro algebra is the ``Euler
number" $|\G^*|$, which is the same as the rank of the lattice
$R_\Z(\G)$.
\end{remark}

\begin{theorem}  \label{th:cubic}
For $\be \in \Clg$, the operator $\mathfrak b (\be)$ is given by
the zero-mode of a vertex operator:
$$\mathfrak b (\be) = \frac16
:\mathfrak p^3:_0 (\tau_{3*}\be).$$
In particular, for $\beta\in
\G^*$, we have $\mathfrak b (\be) = \frac16 h_\beta^2 :\mathfrak p
(\beta)(z)^3:_0.$
\end{theorem}

\begin{proof}
By applying Lemma~\ref{lem:push} twice, we obtain for $\be \in \G^*$ that
 \begin{eqnarray*}
 \tau_{3*} \be
  = (\tau_{2*} \otimes \Id) (h_\be  \be \otimes \be)
  = h_\be^2 \; \be \otimes \be \otimes \be.
  \end{eqnarray*}
It follows that the two identities to be proved are equivalent
thanks to the linearity on $\be$.

Since the operators $\mathfrak b (\be)$ and $h_\be^2/6  :\mathfrak p(\be)(z)^3:_0$
annihilate the vacuum vector $\vac$, it
suffices to see that they have the same commutators with
$\mathfrak b_n(\g), \g\in \G^*$.

We calculate that
\begin{eqnarray*}
 &&\left[ h_\be^2/6  :\mathfrak p(\be)(z)^3:_0, \mathfrak p_n (\g) \right] \\
 &=&  h_\be^2 \left[1/2\cdot \sum_{k>0,m>0} (\mathfrak p_{-k}(\be)
   \mathfrak p_{-m}(\be) \mathfrak p_{k+m}(\be)
   + \mathfrak p_{-k-m}(\be) \mathfrak p_k(\be) \mathfrak p_{m}(\be)), \mathfrak
   p_n(\g) \right]\\
 &=& h_\be^2 \cdot (-n) \cdot \delta_{\be,\g} L_n^\be
 = -n h_\be \delta_{\be,\g} \, L_n(\be)
 = -n \, L_n(\be  \g).
\end{eqnarray*}
This coincides with the commutator $[\mathfrak b(\be), \mathfrak
p(\g)]$ as given in Theorem~\ref{th:virasoro}.
\end{proof}

\begin{remark}
Theorems~\ref{th:virasoro} and \ref{th:cubic} are reformulations
of Theorems 4 and 3 in \cite{FW}.
In the symmetric group case, i.e. $\G$ is trivial, the operator
$\mathfrak b$ was first considered by Goulden \cite{Gou} for a
different purpose and he established Theorem~\ref{th:cubic} (in a
different form) in the setup of ring of symmetric functions.
This operator was rediscovered in
\cite{FW} in the general wreath product setup
for the group theoretic construction of the Virasoro
algebra. In our new formulation above, we have made the transfer
property of these commutators transparent in a way
parallel to the developments in the Hilbert schemes and symmetric
products (cf. \cite{QW} and the references therein).
\end{remark}

\subsection{The general convolution operators and Heisenberg algebra}

Given $\mathfrak f \in \End (\RG)$, we denote
$\mathfrak f' =[\mathfrak b, \mathfrak f] ={\rm ad} \,\mathfrak b (\mathfrak f)$
and define inductively
$\mathfrak f^{(k)} = {\rm ad} \,\mathfrak b (\mathfrak f^{(k-1)}).$

\begin{theorem} \label{th:covcomm}
Let $\g, \alpha \in \Clg$. Then we have
\begin{eqnarray*}
 [\mathfrak O_{\hbar}(\g), \mathfrak p_{-1}(\alpha)]
 &=&\exp(\hbar \cdot {\rm ad}\,{\mathfrak b}) (\mathfrak p_{-1}(\g \alpha))
 \\
 {[} \mathfrak O^k(\g), \mathfrak p_{-1}(\alpha)] &=& \mathfrak p^{(k)}_{-1}(\g\alpha),\quad k \ge 0.
\end{eqnarray*}
\end{theorem}

\begin{proof}
The two identities are equivalent, and we will prove the second identity.

Recall that $\mathfrak p_{-1}(\alpha)(y) =\frac1{(n-1)!} \sum_{g
\in S_n} \ad g\,(y \otimes \alpha)$, for $y \in R(\G_{n-1}).$
Regarding $\G_{n-1}$ as a subgroup of $\Gn$, we have an injective
algebra homomorphism $\iota :\C[\G_{n-1}] \rightarrow \C[\G_n]$
given by the natural inclusion. Thus, using (\ref{eq:average})
we calculate that

\begin{eqnarray*}
 && (n-1)! \;{[} \mathfrak O^k(\g), \mathfrak p_{-1}(\alpha)] (y)  \\
 &=& (n-1)! \;[ \mathfrak O^k(\g)\cdot \mathfrak p_{-1}(\alpha)(y)
  - \mathfrak p_{-1}(\alpha)\cdot \mathfrak O^k(\g)(y)]  \\
 &=& \Xi_n^k(\g) \circ \sum_{g \in S_n} \ad g \,(y \otimes \alpha)
   -\sum_{g \in S_n} \ad g \,((\Xi_{n-1}^k(\g) \circ y) \otimes \alpha)  \\
 &=& \sum_{g \in S_n} \ad g \, [(\Xi_n^k(\g) -\iota(\Xi_{n-1}^k(\g)))\circ (y \otimes \alpha)]
\end{eqnarray*}
where we used the fact that $\Xi_n^k(\g)$ is $S_n$-invariant.
We have $\Xi_n^k(\g) -\iota(\Xi_{n-1}^k(\g)) =
 \xi_{n;n}^{k} \circ \g^{(n)}$, by definition. Thus, we obtain that
\begin{eqnarray*}
 (n-1)! \;{[} \mathfrak O^k(\g), \mathfrak p_{-1}(\alpha)] (y)
 &=& \sum_g \ad g \,[ \xi_{n;n}^{k} \circ \g^{(n)} \circ (y \otimes \alpha)]  \\
 &=& \sum_g \ad g \,[ \xi_{n;n}^{k}  \circ (y \otimes \g\alpha)].
 \end{eqnarray*}

It remains to prove that
\begin{eqnarray} \label{orb:eqinduction}
\sum_{g\in S_n} \ad g \,[ \xi_{n;n}^{k}  \circ (y \otimes
\g\alpha)] =(n-1)! \;\mathfrak p^{(k)}_{-1}(\g\alpha) (y).
\end{eqnarray}
We will proceed by induction. It is trivial for $k=0$.
Note that $\Xi_n^1 -\iota( \Xi_{n-1}^1) = \xi_{n;n}.$ Under
the assumption that the formula (\ref{orb:eqinduction}) is true
for $k$, we have
 \begin{eqnarray*}
 &&  \sum_{g\in S_n} \ad g \,[ \xi_{n;n}^{(k+1)} \circ (y \otimes \g\alpha)]\\
 &=&\sum_g \ad g \,[(\Xi_n^1 -\iota( \Xi_{n-1}^1))
    \circ ( \xi_{n;n})^{ k} \circ (y \otimes \g \alpha)]  \\
 &=&\Xi_n^1 \circ\sum_{g\in S_n} \ad g \,[ \xi_{n;n}^{k} \circ
   (y \otimes \g \alpha)]
   -\sum_{g\in S_n} \ad g \,[\iota (\Xi_{n-1}^1) \circ \xi_{n;n}^{k}
   \circ (y \otimes \g \alpha)],
 \end{eqnarray*}
since $\Xi_n^1(\g)$ is $S_n$-invariant. By using the induction
assumption twice, we get
  \begin{eqnarray*}
  && \sum_{g\in S_n} \ad g \,[ \xi_{n;n}^{(k+1)} \circ (y \otimes \g\alpha)]\\
  &=& (n-1)! \; \Xi_n^1 \circ \mathfrak
  p^{(k)}_{-1}(\g\alpha)(y)
   -\sum_{g\in S_n} \ad g \,[ \xi_{n;n}^{k} \circ ( (\Xi_{n-1}^1 \circ y) \otimes\g \alpha)] \\
  &=& (n-1)!\; [\mathfrak b \cdot \mathfrak p^{(k)}_{-1}(\g\alpha) (y)
   -\mathfrak p^{(k)}_{-1}(\g\alpha) (\Xi_{n-1}^1 \circ y) ] \\
   &=& (n-1)! \;\mathfrak p^{(k+1)}_{-1}(\g\alpha) (y) .
 \end{eqnarray*}
So by induction, we have established (\ref{orb:eqinduction}) and
thus the theorem.
\end{proof}
\section{The $\Winf$ algebra and the Fock space $\RG$}
\label{sec:winf}
\subsection{The $\Winf$ algebra parametrized by $\G^*$}

Let $A$ be a finite-dimensional complex vector space endowed with
a commutative algebra structure with identity. We further assume that
there exists a linear operator $\Tr :A \rightarrow \C$.

Let $t$ be an indeterminate and let $\partial_t=\frac{d}{dt}$. We
denote by ${\D}_{as}$ the associative algebra of regular
differential operators on a circle. Denote by ${\DA}_{as}$ the
associative algebra ${\D}_{as} \otimes A$. If $a_1, \ldots, a_N$
is a linear basis for $A$, then
$t^{\ell+k} (\partial_t)^\ell \otimes a_i$, where $\ell \in \Z_+, k\in\Z$,
$1 \le i\le N$, form a linear basis for ${\DA}_{as}$.

Let $\DA$ denote the Lie algebra obtained from ${\DA}_{as}$ by
taking the usual Lie bracket:

$$[X\otimes \alpha, Y \otimes \beta] = (XY -YX) \otimes (\alpha\beta ).$$

Denote by $\hDA$ the central extension of $\DA$ by a
one-dimensional vector space with a specified generator $C$:
\begin{eqnarray*}
 0 \longrightarrow \C C \longrightarrow \hDA \longrightarrow
 {\DA} \longrightarrow 0.
\end{eqnarray*}
The commutation relation in $\hDA$ is given by
\begin{eqnarray}
 && \left[
     t^r f(D) \otimes \alpha, t^s g(D) \otimes \beta
  \right] \nonumber \\
    & = & t^{r+s}
    \left(
      f(D + s)g(D) - f(D)g(D +r)
    \right) \otimes (\alpha\beta) \nonumber \\
   & + & \Psi
      \left(
        t^r f(D), t^s g(D)
      \right)
      C,
  \label{eq:winfcom}
\end{eqnarray}
where $\Psi: \DA \times \DA \rightarrow \C$ is the $2$-cocycle
given by:

\begin{eqnarray*} \label{eq_cocycle}
  \Psi
      \left(
        t^r f(D) \otimes \alpha, t^s g(D) \otimes \beta
      \right)
   & = &
   \left\{
      \everymath{\displaystyle}
      \begin{array}{ll}
         h_\be^2 \Tr( \alpha \beta) \sum_{-r \leq j \leq -1} f(j) g(j+r)
           , & r = -s \geq 0  \\
          0, & r + s \neq 0.
      \end{array}
    \right.
  \label{cocy}
\end{eqnarray*}

In the case when $A =\Clg$ which we are concerned, we shall write
$\hD(\Clg)$ as $\hDg$. When $\G$ is trivial and $A =\C$, we will
simply write $\hDg$ as $\hD$, which is the usual $\Winf$ algebra,
cf. \cite{FKRW, Kac}. The algebra $\hDg$ has a linear basis given
by
\begin{eqnarray*}
  J^l_k (\alpha)= - t^{l+k} ( \partial_t )^l \otimes h_\alpha^{-1} \alpha, \quad
    l \in \Z_{+}, k \in \Z, \alpha \in \G^*.
\end{eqnarray*}
Set $D = t \partial_t$. A different basis of $\hDg$ is given by
\begin{eqnarray*}
 L^l_k (\alpha) =
 - t^{k} D^l \otimes h_{\alpha}^{-1} \alpha,  \quad l \in \Z_{+}, k \in \Z, \alpha \in \G^*.
\end{eqnarray*}
 Note that $f(D) t = t f( D +1)$ for
$f(w) \in \C [w]$ and hence $J^l_k(\alpha) = -t^k [D]_l \otimes h_\alpha^{-1} \alpha$,
where we have used the notation
$ [x]_l = x(x-1)\ldots (x-l+1).$
We then extend $J^l_k (\alpha)$ and $L^l_k (\alpha)$ by linearity to all
$\alpha \in \Clg$.
We further introduce the vertex operator
\begin{eqnarray*}
 J^l(\alpha)(z) &= &\sum_{k \in\Z}J^l_k(\alpha) z^{-k-l-1} .
\end{eqnarray*}
In particular, the components of $J^0(z)$ generates a Heisenberg
algebra while the components of $J^1(z)$ generates a Virasoro
algebra. The algebra $\hDg$ has a weight gradation given by $wt
(J^l_k(\alpha)) =-k$. The Cartan subalgebra of $\hDg$ has a linear
basis given by the center $C$ and all polynomials in $D$.

\begin{remark} \label{rem:good}
 We have introduced $J^\ell_k (\alpha)$
and the $2$-cocycle $\Psi$ by taking into account the fact that
$h_\alpha^{-1}\g$, where $\g \in \G^*$, are orthogonal idempotents.
For $\G$ trivial, we simply ignore all the parameters lying in
$\Clg$ (or equivalently set them to be 1).
\end{remark}
\subsection{Free field realization}

Take a pair of free fermionic fields
 $$\psi^+ (z) =\sum_{n\in \Z+\hf} \psi^{+}_n z^{-n-\hf},
 \quad \psi^- (z) =\sum_{n\in \Z+\hf} \psi^{-}_{n}z^{-n-\hf}$$
where $\psi^{\pm}_{n}$ satisfied the following anti-commutation
relations:
 $$[\psi^{+}_{m}, \psi^{-}_{n}]_+ =\delta_{m,-n}$$
and the other commutators being zero. We denote by $\FFock$ the
fermionic Fock space generated by the vacuum vector $\vac$ from
the creation operators $\psi^{\pm}_{n}, n <0.$ We have the
(charge) decomposition $\FFock = \sum_{k \in \Z} \FFock^{(k)}$,
where $\FFock^{(k)}$ is spanned by the monomial
$\psi^{+}_{m_1}\psi^{+}_{m_2}\cdots\psi^{+}_{m_p}
\psi^{-}_{n_1}\psi^{-}_{n_2}\cdots\psi^{-}_{n_q} \cdot \vac,$ where $ p-q
=k.$

It is known \cite{Mat, FKRW} that one can realize $\hD$ as
\begin{eqnarray*}
J^l(z) = :\partial_z^l \psi^-(z) \psi^+(z): , \quad l \ge 0.
\end{eqnarray*}
In particular $J^0(z)$ is a Heisenberg vertex operator (i.e. a
free boson). This is part of the classical boson-fermion
correspondence.

It is
a well-known fact that $\FFock^{(k)}$ is irreducible under the
action of the Heisenberg algebra.
As $\hD$ contains the Heisenberg algebra generated
by the field $J^0(z)$, it also acts on $\FFock^{(k)}$ for each $k$
irreducibly. We will be only interested in $\FFock^{(0)}$.

The boson-fermion correspondence further allows us to write
$\psi^{\pm}(z)$ in terms of the free bosonic field $J^0(z)$ (cf.
e.g. \cite{Kac}):

\begin{eqnarray} \label{eq_corres}
 \psi^{\pm}(z)
 &=& : \exp \left(\int \mp J^0(z)dz \right): S_{\pm } \nonumber\\
 &= & \exp (\mp \sum_{k>0} J^0_{-k}z^k/k) \exp(\pm \sum_{k>0}
 J^0_{k}z^{-k}/k) z^{\mp J^0_0} S_{\pm},
\end{eqnarray}
where $S_{\pm}$ is the shifted operator $\FFock^{(k)} \rightarrow
\FFock^{(k \pm 1)}$ which commutes with the action of the
Heisenberg algebra.

\begin{proposition} \label{prop:bosoniz}
On $\FFock^{(0)}$, we can realize the fields $J^{l-1}(z)$ $(l \ge
1$) as  normally ordered polynomials, denoted by $P_l(J^0)$ or
$P_l(J^0(z))$, in terms of $J^0(z)$ and its derivative fields.
More precisely, we have $J^{l-1}(z) =\frac1l P_l(J^0(z))$, where
\begin{eqnarray*}
 P_l(J^0(z)) = \frac{\partial^l : \exp(\int  J^0(z)dz):}{: \exp(\int  J^0(z)dz):}.
\end{eqnarray*}
\end{proposition}

For example, we can write down explicitly for small $l$:
\begin{eqnarray*}
   P_1(J^0)&=&  J^0\\
   P_2(J^0)&=& :(J^0)^2: +\partial J^0\\
   P_3(J^0)&=& :(J^0)^3: +3J^0 \partial J^0 +\partial^2J^0.
\end{eqnarray*}

\begin{proof}
By the boson-fermion correspondence (\ref{eq_corres}) and the
identity
\begin{eqnarray*}
 \left[- \sum_{k>0}J^0_{k}z^{-k}/k,\sum_{k>0} J^0_{-k} w^k/k \right]
 =-\sum_{k>0} \frac1k (w/z)^k = \ln(1-w/z),
\end{eqnarray*}
we have the following standard computation of the
operator product expansion for $|z| >|w|$:

\begin{eqnarray*}
 && \psi^-(z)\psi^+(w) \\
 &=& \frac1{z-w} \exp \left( \sum_{k>0} J^0_{-k}(z^k -w^k)/k \right)
 \exp \left(- \sum_{k>0} J^0_{k}(z^{-k}-w^{-k})/k \right) \\
 &=& \frac1{z-w}  \sum_{l \ge 0}\frac1{l!} (z-w)^l \cdot \\
 &&\quad \cdot \left[\partial_z^l
 \exp ( \sum_{k>0} J^0_{-k}(z^k -w^k)/k) \exp(- \sum_{k>0}
 J^0_{k}(z^{-k}-w^{-k})/k) \right]|_{z=w} \\
 &=& \frac1{z-w} + \sum_{l \ge 1}\frac1{l!}(z-w)^{l-1} P_l(J^0(w)).
\end{eqnarray*}

We can calculate the operator product expansion
$\psi^-(z)\psi^+(w)$ for $|z| >|w|$ in a different way:

\begin{eqnarray*}
\psi^-(z)\psi^+(w)
 &=& \frac1{z-w} + : \psi^-(z)\psi^+(w): \\
 &=& \frac1{z-w} + \sum_{k \ge 0} \frac1{k!}  (z-w)^k :\partial_w^k
 \psi^-(w)\psi^+(w):\\
 &=& \frac1{z-w} + \sum_{k \ge 0} \frac1{k!}  (z-w)^k J^k(w).
\end{eqnarray*}

The proposition follows by comparing the right-hand sides of the
above two calculations.
\end{proof}

 For each $\g\in \G^*$, we introduce a pair of
fermionic fields
$$\psi^{\pm}(\g)(z)
=\sum_{n\in \Z+\hf} \psi^{\pm}_n(\g) z^{-n-\hf},$$
where $\psi^{\pm}_{n}(\g)$ satisfies the following anti-commutation
relations:
 $$[\psi^{+}_{m}(\be), \psi^{-}_{n}(\g)]_+ =\delta_{m,-n} \delta_{\be,\g}, \quad \be,\g \in \G^*,$$
and the other commutators being zero. We denote by $\FFock_\G$ the
corresponding fermionic Fock space associated to $ \psi^{\rm}(\g), \g \in \G^*$,
and similarly we have the charge decomposition
$\FFock_\G = \oplus_{k\in \Z} \FFock_\G^{(k)}$.
Now, by Remark~\ref{rem:good}, we can realize $\Winf^\G$ acting on $\FFock_\G$ as
\begin{eqnarray*}
J^l(\g)(z) = :\partial_z^l \psi^-(\g)(z) \psi^+(\g) (z): , \quad \g \in \G^*.
\end{eqnarray*}
Of course, Proposition~\ref{prop:bosoniz} applies when we replace
$J^l(z)$ by $J^l(\g)(z)$.
\subsection{Convolution operators and $\hDg$}

Note that the Heisenberg algebra generated by $J^0_n(\g)$
and also $\hDg$ acts irreducibly on $\FFock_\G^{(0)}$.
Since the components of $J^0(\g)(z)$ and those of $\mathfrak p(\g) (z)$ satisfy
the same Heisenberg algebra commutation relation, we may
identify $\RG$ with $\FFock_\G^{(0)}$ and let $J^0_n(\g)$
act on $\RG$ as $\mathfrak p_n(\g)$. The
point here is how to realize the action of $\hDg$ on $\RG$ in a
group-theoretic manner.

\begin{proposition}  \label{prop:quardiff}
With the above identifications, the operator $\mathfrak b(\g)$ for $\g
\in \G^*$ can be identified with $- h_\g \cdot
D (D+1)/2 \otimes \g$ in $\hDg.$
\end{proposition}

\begin{proof}
Recall that $J^0_0(\g) =\mathfrak p_0(\g)$ acts on $\RG$ as 0.
  From Theorem~\ref{th:cubic} and Proposition~\ref{prop:bosoniz}, we see that
 \begin{eqnarray*}
 \mathfrak b(\g)
 &=& \frac16 h_\g^2 :\mathfrak p (\g)^3:_0 \\
 &=&  h_\g^2 \cdot \left(P_3(J^0 (\g))/6 - \pa P_2(J^0(\g))/4 \right)  \\
 &=&  h_\g^2 \cdot(J^2_0(\g)/2 - \pa J^1_0(\g)/2) \\
  &=& h_\g^2 \cdot(-D (D -1)/2 - D) \otimes h_\g^{-1} \g \\
 &=& - h_\g \cdot D (D+1)/2 \otimes \g.
  \end{eqnarray*}
\end{proof}

In particular, in the symmetric group case, we have
$\mathfrak b =-D (D+1)/2.$

\begin{proposition}  \label{prop:convdiff}
Let $q=e^\hbar$. We can identify the operator ${\mathfrak
O}_{\hbar}(\g)$, where $\g\in \G^*$, on $\RG$ as the following
weight-zero differential operator in $\hDg$:
 $$ {\mathfrak O}_{\hbar} (\g)
 = \frac{q^{h_\g^2\, D}-1}{q^{-h_\g^2}-1} \otimes h_\g^{-1}\g.$$
\end{proposition}

\begin{proof}
 From (\ref{eq:winfcom}) and Proposition~\ref{prop:quardiff}, for $\g \in \G^*$, we have
 $$  [\mathfrak b(\g), \mathfrak p_{-1}(\be)] =
 [-h_\g D (D+1)/2 \otimes h_\g^{-1}\g, t^{-1} \otimes h_\be^{-1}\be]
 = t^{-1}( h_\g^2D) \otimes h_\g^{-1}h_\be^{-1}\g\be.$$
 Together with Theorem~\ref{th:covcomm}, this implies that (for $\be\in \G^*$)

 $$
[{\mathfrak O}_{\hbar}(\g), \mathfrak p_{-1}(\be)]
 = q^{ \mbox{ad}\, \mathfrak b(\g)} \mathfrak p_{-1}(\be)
 = q^{ \mbox{ad}\, \mathfrak b(\g)} (t^{-1} \otimes h_\be^{-1} \be)
 = t^{-1}q^{h_\g^2 D} \otimes (h_\g^{-1} h_\be^{-1}\g\be),$$
where we have used the fact that $h_\g^{-1}\g$ are orthogonal
idempotents and $$(h_\g^{-1}\g)^n (h_\be^{-1} \be) =h_\g^{-1} h_\be^{-1}\g\be, \quad n \ge 1. $$

 On the other hand, a simple calculation using (\ref{eq:winfcom})
 gives us
 \begin{eqnarray*}
 && \left[(q^{-h_\g^2}-1)^{-1}(q^{h_\g^2 D}-1) \otimes h_\g^{-1}\g,
 \mathfrak p_{-1}(\be) \right] \\
 &=& (q^{-h_\g^2}-1)^{-1} t^{-1} (q^{h_\g^2 (D-1)} -q^{h_\g^2 D})
 \otimes (h_\g^{-1} h_\be^{-1} \g\be)
 = t^{-1}q^{h_\g^2 D} \otimes (h_\g^{-1} h_\be^{-1} \g\be).
\end{eqnarray*}

 Note that both operators ${\mathfrak O}_{\hbar}(\g)$ and
$(q^{-h_\g^2}-1)^{-1}(q^{h_\g^2 D}-1) \otimes h_\g^{-1}\g$
 commutes with $\mathfrak b (\be)$, and both operators
 annihilates the vacuum vector. It follows from the following simple
Lemma~\ref{lem:triv} that these two operators coincides.
\end{proof}

\begin{corollary}
In the case of symmetric groups (i.e. $\G$ is trivial),
we identify the operator ${\mathfrak O}_{\hbar}$ as the
degree-zero differential operator $\frac{q^D-1}{q^{-1}-1}$ in $\hD$.
\end{corollary}

\begin{lemma} \label{lem:triv}
Given two operators $\mathfrak f_1, \mathfrak f_2 \in \End(\RG)$
such that $\mathfrak f_1 \vac =\mathfrak f_2 \vac=0$ and
$[\mathfrak f_1, \mathfrak b] =[\mathfrak f_2, \mathfrak
b]=0,$ and $[\mathfrak f_1, \mathfrak p_{-1}(\be)] =[\mathfrak
f_2, \mathfrak p_{-1}(\be)]$ for all $\be \in \Clg$. Then
$\mathfrak f_1=\mathfrak f_2$.
\end{lemma}

\begin{proof}
Follows from the fact that we can obtain the whole $\RG$ by repeatedly
applying the operator $\mathfrak b$ and $\mathfrak p_{-1}(\g)$,
$\g \in \Clg$.
\end{proof}

\begin{theorem}  \label{th:winf}
The commutators $[ \mathfrak O^k (1), \mathfrak p_n(\alpha)]$, where
$k \ge 0, n \in\Z,$ and $ \alpha \in \G^*$,
realize a level one
irreducible representation of the Lie algebra $\hDg$ on $\RG$.
\end{theorem}

\begin{proof}
If we expand the formula for $\mathfrak O_\hbar (\g)$
in Proposition~\ref{prop:convdiff} as a Taylor series in $\hbar$, we see
that $\mathfrak O^k(\g)$ is of degree $k+1$ in
 $D$.
Identifying $\mathfrak p_n(\alpha)$ with $t^{-n} \otimes h_\alpha^{-1}\alpha$,
we see that the commutators $[ \mathfrak O^k (1), \mathfrak p_n(\alpha)]$
is $t^{-n}$ times a polynomial in $D$ of degree $k$ and they form
a linear basis for $\hDg$. It is clear that the representation is of level one.
\end{proof}

\subsection{Convolution operators and vertex operators}

\begin{lemma}  \label{lem_variable}
 Let $q=e^{\hbar}$. We have:
 \begin{eqnarray*}
 q^D &=&1+  \sum_{l \ge 1} \frac1{l!}  (q-1)^l {[D]}_l
 = 1+  \sum_{l \ge 1} \frac1{l!} \hbar^l D^l, \\
 \frac{q^D-1}{q^{-1}-1} &=& -q\sum_{l \ge 1} \frac1{l!} (q-1)^{l-1} [D]_l .
 \end{eqnarray*}
\end{lemma}

\begin{proof}
Follows from the expansions of $q^D$ and $(q^D-1)/(q^{-1}-1)$ in
terms of powers of $(q-1), (q^{-1}-1),$ and $\hbar$ respectively.


\end{proof}

Let us introduce the following vertex operator associated to $\g\in \G^*$:
\begin{eqnarray*}
 V(\g;z,q) &=& (q-1)z\,\psi^-(\g) (qz) \psi^+(\g)(z) \\
 &=& \exp \left(\sum_{k > 0} \frac{(q^{k}-1)z^{k}}{k} \mathfrak p_{-k}(\g)
          \right)
 \exp \left(\sum_{k > 0} \frac{(1-q^{-k})z^{-k}}{k} \mathfrak p_k(\g)
      \right).
\end{eqnarray*}
We write $V(\g;z,q) = \sum_{m \in \Z} V_m(\g;q)z^{-m}$.

The identification of $\mathfrak O_\hbar (\g)$ as an differential operator
is useful as illustrated in the proof of the following theorem which
describes the connection between JM elements and vertex operator.

\begin{theorem} \label{th:vo}
 Let $q=e^{\hbar}$, and $\g \in \G^*.$ The operator $\mathfrak O_\hbar$ can be
 expressed as:
 \begin{eqnarray*}
 {\mathfrak O}_\hbar(\g)  =\frac{q}{(q-1)^2} \left( V_0(\g; q^{h_\g^2})-1 \right).
 \end{eqnarray*}
\end{theorem}

\begin{proof}
As usual, we identify $J^0_k(\g)$ with $ \mathfrak p_k(\g)$.
By taking the Taylor expansion of $V(\g;z,q)$ with respect to
$(q-1)$ and using Proposition~\ref{prop:bosoniz}, we obtain:

\begin{eqnarray} \label{eq_taylor}
 V(\g;z,q)
 &=& 1+ \sum_{k \ge 1} \frac1{k!} (q-1)^k P_k(J^0(\g))z^k \nonumber \\
 &=& 1+ \sum_{l \ge 0} \frac1{l!} (q-1)^{l+1} J^{l}(\g)(z) z^{l+1}.
\end{eqnarray}

Note that the operator $J^0_0(\g)$ acts as $0$,
and $J^l_0(\g) =-[D]_l \otimes h_\g^{-1} \g$. By (\ref{eq_taylor}) and
Lemma~\ref{lem_variable} (where $q$ is replaced by $q^{h_\g^2}$), we have
\begin{eqnarray*}
  \frac{q^{h_\g^2}}{(q^{h_\g^2}-1)^2} \left(V_0(\g;q^{h_\g^2})-1 \right)
 &=& \frac{q^{h_\g^2}}{(q^{h_\g^2}-1)^2}
   \sum_{l \ge 1} \frac1{l!} (q^{h_\g^2}-1)^{l+1} \cdot J^{l}_0(\g)  \\
 &=& - q^{h_\g^2}
 \sum_{l \ge 1} \frac1{l!} (q^{h_\g^2}-1)^{l-1}\cdot [D]_l \otimes h_\g^{-1} \g  \\
 &=&  \frac{q^{h_\g^2 \cdot D}-1}{q^{-h_\g^2}-1} \otimes h_\g^{-1}\g.
\end{eqnarray*}
Now our result follows from comparing with
Proposition~\ref{prop:convdiff}.
\end{proof}

\begin{remark}
In the symmetric group case, Theorem~\ref{th:winf} and Theorem~\ref{th:vo} have been
established by Lascoux and Thibon in a
different approach \cite{LT} (as a generalization of the construction of the Virasoro
algebra by I.~Frenkel and the author \cite{FW}).
Actually, we have been unsuccessfully seeking for the $\Winf$ algebra
acting on $\RG$ right after we found the Virasoro algebra construction.
In particular, we noticed the relevance of the vertex operator $V(z,q)$, but did not
come up with the right group-theoretic construction (i.e. the JM elements) at the time.
\end{remark}
\section{Applications: stability and algebra generators for $\Rgn$}
\label{sec:application}

An observation made in \cite{QW} is the following. Assume that ({\bf A1}) there
exists a sequence of (finite-dimensional) Frobenius
$\C$-algebras $A^{[n]}$ ($n \ge 0$) such that $A^{[0]} =\C$ and $A=A^{[1]}$. ({\bf A2})
the direct sum $\oplus_n A^{[n]}$ affords the structure of a Fock
space of a Heisenberg algebra as in Theorem~\ref{th:heis}
(where $\Clg$ is replaced by $A$).
({\bf A3}) There exists a sequence
of elements $\Xi_n^k(\alpha) \in A^{[n]}$ depending on $\alpha \in
A$ (linearly) and a non-negative integer $k$, which can be used to
define operators $\mathfrak O^k(\alpha)$. The operators $\mathfrak
O^k(\alpha)$, $\mathfrak b (\alpha) =\mathfrak O^1(\alpha)$ and the Heisenberg algebra generators
satisfy the relations as in Theorem~\ref{th:cubic} and Theorem~\ref{th:covcomm}
(where $\Clg$ is replaced by $A$).
Then various statements, such as the stability and algebra generators,
can be derived from these three axioms by a standard procedure
using only the Heisenberg
generators, the diagonal pushforward maps, the Euler class, etc.
This was first done for cohomology rings of Hilbert schemes of points
on projective surfaces \cite{LQW}, and then it has been applied to orbifold
cohomology rings of the symmetric products \cite{QW}.

In this section, we will formulate some consequences of these three
axioms, which are the counterparts of the results obtained
in \cite{LQW} and in Section~4, \cite{QW}.
We will mainly formulate the stability result
which allow us to introduce the {\em stable algebra}
$\mathfrak A_\G$ associated
to $\G$. The proofs are the same as in
\cite{LQW} (also cf. \cite{QW}) and thus will not be repeated here.
There are several more consequences of these axioms
which we do not reproduce here, cf. \cite{LQW}.
We further follow the ideas of Ivanov and Kerov \cite{IK}
to establish that the structure constants of $\mathfrak A_\G$
are nonnegative integers.
\subsection{The stability in class algebras $\Rgn$}

Let $s \ge 1$, and let $\alpha_1, \ldots, \alpha_s \in \Clg$.
Let $\pi = \{\pi_1, \ldots, \pi_j \}$ be a partition of
the set $\{1, \ldots, s\}$, and define $\ell(\pi)=j$,
$\alpha_{\pi_i} =\prod_{m\in \pi_i} \alpha_m$.
We denote by ${\bf
1}_{-k} = \frac{{\mathfrak p}_{-1}(1)^k}{k!}$ if $k \ge 0$ and
${\bf 1}_{-k} = 0$ if $k<0.$

\begin{proposition}
Let $n, s \ge 1$, $k_1, \ldots, k_s \ge 0$, and let $\alpha_1,
\ldots, \alpha_s \in \Clg$. Then, the convolution
product $ \Xi_n^{k_1}(\alpha_1) \circ \cdots \circ
\Xi_n^{k_s}(\alpha_s)$ in $\Rgn$ is a finite linear combination of
expressions of the form

\begin{eqnarray*}
{\bf 1}_{-\left (n-\sum\limits_{i=1}^{\ell(\pi)}
  \sum\limits_{j=1}^{m_i-2r_i} n_{i, j} \right )}
\prod_{i = 1}^{\ell(\pi)} \left ( \prod_{j = 1}^{m_i-2r_i}
\mathfrak p_{-n_{i, j}} \right ) (\tau_{(m_i-2r_i)*}(\chi^{r_i}
\alpha_{\pi_i})) \cdot |0\rangle
\end{eqnarray*}
where $\pi$ runs over all partitions of $\{1, \ldots, s \}$,
$m_i, r_i \in \Z_+$ such that
$$2r_i \le m_i \le 2+ \sum_{j \in \pi_i} k_j,$$
$0 < n_{i, 1} \le \ldots \le n_{i, m_i-2r_i}$,
$\sum\limits_{j=1}^{m_i-2r_i} n_{i, j} \le \sum\limits_{j \in
\pi_i} (k_j+1)$ for every $i$, and

\begin{eqnarray*}
\sum_{i = 1}^{\ell(\pi)} \left ( m_i - 2 + \sum_{j = 1}^{m_i-2r_i}
n_{i, j} \right ) = \sum_{i=1}^s k_i.
\end{eqnarray*}
Moreover, all the coefficients in this linear combination are
independent of the group $\G$, $\alpha_1,
\ldots, \alpha_s$, and the integer $n$.
\end{proposition}

\begin{remark}
This proposition is the counterpart of a theorem in \cite{LQW}.
\end{remark}

\begin{theorem} \label{th:cupprod}
Let $s \ge 1$ and $k_i \ge 1$ for $1 \le i \le s$. Fix $n_{i, j}
\ge 1$ and $\alpha_{i, j} \in \Clg$ for $1 \le j \le k_i$, and let
$n \ge \sum\limits_{j=1}^{k_i} n_{i, j}$ for all $1 \le i \le s$.
Then the convolution product
\begin{eqnarray*}
\prod_{i=1}^s \left ( {\bf 1}_{-(n - \sum_{j=1}^{k_i} n_{i, j})}
\prod_{j=1}^{k_i} \mathfrak p_{-n_{i, j}}(\alpha_{i, j}) \cdot
|0\rangle \right )
\end{eqnarray*}
in $\Rgn$ is equal to a finite linear combination of monomials of
the form
\begin{eqnarray*}
 {\bf 1}_{-(n - \sum_{a=1}^N m_{a})}
 \prod_{a=1}^N \mathfrak p_{-m_{a}}(\g_{a}) \cdot |0\rangle
\end{eqnarray*}
where $\sum_{a=1}^N m_a \le \sum_{i=1}^s \sum_{j=1}^{k_i}
n_{i,j}$, and $\g_1, \ldots, \g_N$ depend only on $\chi,
\alpha_{i,j}$, $1\le i \le s, 1\le j\le k_i$. Moreover, the
coefficients in this linear combination are independent of $\G$,
$\alpha_{i,j}$ and $n$.
\end{theorem}
\subsection{The stable algebra $\mathfrak A_\G$}



For a given
$\rho =(\rho (c))_{c \in \G_*} \in {\mathcal P}(\G_*)$,
where $\rho(c) =(r^{m_r(c)})_{r \ge 1} =(1^{m_1(c)} 2^{m_2(c)}
\ldots )$, we define

\begin{eqnarray*}
  {\mathfrak p}_{- \rho(c)}(K^c)
  &=& \prod_{r \ge 1}
  {\mathfrak p}_{-r}(c)^{m_r(c)} = {\mathfrak p}_{-1}(K^c)^{m_1(c)}
  {\mathfrak p}_{-2}(K^c)^{m_2(c)} \cdots    \\
 {\mathfrak p}_{\rho}(n)
 &=& {\bf 1}_{-(n-\Vert \rho \Vert)}
  \prod_{c\in \G_*}{\mathfrak p}_{-{\rho}(c)}(K^c)\cdot\vac \in
  \Rgn, \quad \text{for } n \geq \Vert \rho \Vert
\end{eqnarray*}
and set ${\mathfrak p}_{\rho}(n)=0$ for $n < \Vert \rho \Vert$.
We also define $\widetilde{\mathfrak p}_{\rho}(n)
={\widetilde{z}_\rho }^{-1} {\mathfrak p}_{\rho}(n)$,
where the constant
$\widetilde{z}_\rho =\prod_{r\geq 1,c \in \G_*} r^{m_r(c)} m_r(c)!$
is independent of $n$.

For a given $\rho$ with $\Vert \rho \Vert \leq n$, we denote by
$\widetilde{\rho} = \rho \cup (1^{n - \Vert \rho \Vert}) \in
\mathcal P_n(\G_*)$, where $\widetilde{\rho} (c) ={\rho}(c)$
for $c \neq c^0$ and $\widetilde{\rho} (c^0)$ is obtained from $\rho(c^0)$
by increasing the multiplicity of 1 by $n - \Vert \rho \Vert.$
The Frobenius characteristic map (cf. \cite{Mac, FJW}) allows one to
identify the class function
$ {n - \Vert \rho \Vert +m_1(c^0)\choose m_1(c^0)}
K^{\widetilde{\rho}}$ and ${\mathfrak p}_{\rho}(n)$.

As $\rho$ runs over all partition-valued functions on $\G_*$ with
$\Vert \rho \Vert \le n$, the elements ${\mathfrak
p}_{\rho}(n)$ linearly span $\Rgn$, as a corollary to
Theorem~\ref{th:heis}. According to Theorem~\ref{th:cupprod} (for
$s=2$), we can write the product in the class algebra $\Rgn$ as

\begin{eqnarray} \label{eq_structure}
{\mathfrak p}_{\rho}(n) \circ {\mathfrak p}_{\sigma}(n) =
\sum_{\nu}d_{\rho\sigma}^\nu {\mathfrak p}_{\nu}(n),
\end{eqnarray}
where $\Vert \nu \Vert \le \Vert \rho \Vert +\Vert \sigma \Vert$,
and the structure coefficients $d_{\rho\sigma}^\nu$ are
independent of $n$. Even though the elements ${\mathfrak
p}_{\nu}(n)$ with $\Vert \nu \Vert \le n$ in $\Rgn$ are not
linearly independent, we can show that (cf. \cite{LQW})
the constants $d_{\rho\sigma}^\nu$ in the formula
(\ref{eq_structure}) are uniquely determined from the fact that
they are independent of $n$. Similarly, we can write
uniquely

\begin{eqnarray*} \label{eq_modif}
\widetilde{\mathfrak p}_{\rho}(n) \circ \widetilde{\mathfrak p}_{\sigma}(n) =
\sum_{\nu}\widetilde{d}_{\rho\sigma}^\nu\widetilde {\mathfrak p}_{\nu}(n)
\end{eqnarray*}
where $\Vert \nu \Vert \le \Vert \rho \Vert +\Vert \sigma \Vert$,
and the new structure coefficients $\widetilde{d}_{\rho\sigma}^\nu$ are
also independent of $n$. Indeed, we have
$\widetilde{d}_{\rho\sigma}^\nu =
\widetilde{z}_\nu \widetilde{z}_\rho^{-1}
 \widetilde{z}_\sigma^{-1} d_{\rho\sigma}^\nu. $

The {\it stable algebra} associated to a finite group $\G$, denoted
by ${\mathfrak A}_\G$, is defined to be the algebra with a linear
basis formed by the symbols ${\mathfrak p}_\rho$, $\rho \in
{\mathcal P}(\G_*)$ and with the multiplication defined by
\begin{eqnarray*}
{\mathfrak p}_{\rho} \, {\mathfrak p}_{\sigma} = \sum_{\nu}
d_{\rho\sigma}^\nu {\mathfrak p}_{\nu}
\end{eqnarray*}
where the structure constants $d_{\rho\sigma}^\nu$ are from the
relations (\ref{eq_structure}) and
$\Vert \nu \Vert \le \Vert \rho \Vert +\Vert \sigma \Vert$.

Clearly the stable algebra ${\mathfrak A}_\G$ itself is
commutative and associative. The algebra ${\mathfrak A}_\G$
captures all the information of the class algebra $\Rgn$ for
all $n$, as we easily recover the relations (\ref{eq_structure})
from the algebra ${\mathfrak A}_\G$. We summarize these observations
into the following.

\begin{theorem} {\rm (Stability)} \label{th:stab}
For a finite group $\G$, the product in
the class algebras $\Rgn$ $(n\ge 1)$ can be written in
the form (\ref{eq_structure}) whose structure constants
are independent of $n$. This
give rise to the stable algebra ${\mathfrak A}_\G$ which completely
encodes the class algebra structure of $\Rgn$ for each $n$.
\end{theorem}

\begin{remark}
 In the formulation of the stable algebra, we have
the freedom of choosing a linear basis of $R(\G)$ (cf.
Theorem~\ref{th:cupprod}). In the above formulation,
we have choosen $K^c$ (where $c \in \G_*$) as
the linear basis of $R(\G)$, which allow us to establish
below remarkable integrality and positivity
properties of the structure constants.
\end{remark}

\begin{remark}
Let $\G$ be trivial and thus $\Gn$ reduces to the
symmetric group $S_n$. In this case,
The stability of the class algebras
was first given in Kerov and Olshanski (cf. \cite{KO},
Proposition~3) from a totally different approach
(cf. Lascoux-Thibon \cite{LT}, Section~4, for another proof).
\end{remark}

\begin{remark}
For $\rho \in \mathcal P(\G_*)$ with
$\ell(\rho) = 1$, that is, when the partition $\rho(c)$ is a
one-part partition $(r)$ for some $c \in \G_*$ and is empty
for all the other class functions in $\G_*$, we will simply write $\mathfrak
p_\rho = \mathfrak p_{r,c}$. Following the techniques as developed
in \cite{LQW} (for cohomology rings of Hilbert schemes), we can
show that the stable algebra $\mathfrak A_\G$ is isomorphic to the
polynomial algebra generated by $\mathfrak p_{r,c},\;c\in \G_*,
r \ge 1$.
\end{remark}
\subsection{Positivity and integrality of
the structure constants of $\mathfrak A_\G$}

\begin{theorem}  \label{th:interal}
 The structure constants $\widetilde{d}_{\rho\sigma}^\nu$
for the stable algebra ${\mathfrak A}_\G$
are non-negative integers.
So are $d_{\rho\sigma}^\nu$.
\end{theorem}

\begin{remark}
In the symmetric group case (i.e. $\G$ is trivial),
the integrality and positivity of these structure constants
have been established by Ivanov and Kerov \cite{IK}.
Theorem~\ref{th:interal} also follows from a straightforward
generalization of their ingenious constructions, which we outline below.
In fact, this approach also provides a
second proof of Theorem~\ref{th:stab}.
\end{remark}

All the constructions of Ivanov and Kerov admit a straightforward
generalization to the wreath product setup. We sketch below only
the portion which will provide us
a proof of Theorem~\ref{th:interal}, and refer to \cite{IK} for more detail.

Following \cite{IK}, we introduce the semigroup
of ``partial permutations'' $P\G_n$ as follows.
If $Y$ is a finite set, we denote by $S_Y$ the
symmetric group of permutations on $Y$, and $\G_Y$
the corresponding wreath product. We denote by
$\underline{n}$ the set $\{1,2, \ldots,n\}$,
so $\G_{\underline{n}}$ is just our usual $\Gn$.
For $Y \subset \underline{n}$, we regard $\G_Y$
as a subgroup of $\Gn$.

A {\em partial permutation} of the set $\underline{n}$
is a pair $(Y, a)$ which consists of a finite
subset $Y \subset \underline{n}$ and an element
$a \in \G_Y$. Denote by $P\G_n$ the set of
all partial permutations of $\underline{n}$.
The set $P\G_n$ is endowed with a natural semigroup
structure by letting the
product of two elements $(Y_1, a_1)$ and
$(Y_2, a_2)$ in $P\G_n$ to be
$(Y_1 \cup Y_2, a_1 a_2)$. We denote by
$\Z[P\Gn]$ the semigroup algebra (over $\Z$).

The wreath product $\Gn$ acts on the semigroup $P\G_n$
by $(Y,a) \mapsto (\sigma Y, xax^{-1})$, where $x =(g,\sigma) \in \Gn$,
$g \in \G^n, \sigma \in S_n$. One easily shows that the corresponding
conjugate classes $C_\rho (n)$ in $P\G_n$ are parametrized
by $\rho \in \mathcal P(\G_*)$ such that $\Vert \rho\Vert \leq n$.
By abuse of notation, we will also use $C_\rho (n)$ to denote
the characteristic function of the conjugate class $C_\rho (n)$.
Denote by $R_\Z (P\Gn)$ the $\Z$-span of the characteristic
functions of the conjugate classes in $P\G_n$.

There is a surjective ``forgetful'' ring homomorphism
$\Z[ P\Gn] \rightarrow \Z [\Gn]$ by $(Y,a) \mapsto a$,
and it induces  a surjective ring homomorphism
$R_\Z (P\Gn) \rightarrow R_\Z (\Gn)$ by sending
$C_\rho (n)$ to
${n - \Vert \rho \Vert +m_1(c^0)\choose m_1(c^0)} K^{\widetilde{\rho}}
=\widetilde{\mathfrak p}_\rho (n)$.
It follows that
$$C_\rho (n) \; C_\sigma (n) =
\sum_{\nu} \widetilde{d}_{\rho\sigma}^\nu C_\nu (n)$$
where $\widetilde{d}_{\rho\sigma}^\nu$ are the
structure constants introduced earlier.
This interpretation of $\widetilde{d}_{\rho\sigma}^\nu$
as the structure constants of $R_\Z (P\Gn) $
implies immediately that all $\widetilde{d}_{\rho\sigma}^\nu$
are nonnegative integers.

We can also establish the fact that
$\widetilde{d}_{\rho\sigma}^\nu$ (denoted by
$\widetilde{d}_{\rho\sigma}^\nu(n)$ for the time being) are
independent of $n$ along this line as follows. Let $m \leq n$. We
define a linear map $\theta_{n,m}: \Z[P\Gn] \rightarrow \Z[P\G_m]$
by letting $\theta_{n,m} (Y,a) = (Y,a)$ if $Y \subset
\underline{m}$, and $\theta_{n,m} (Y,a) =0$ otherwise. As in
\cite{IK}, we see that $\theta_{n,m}$ is a surjective ring
homomorphism. Since the homomorphisms $\theta_{n,m}$ for all $m
\leq n$ are compatible, we can define a projective limit of the
algebras $\Z[P\Gn]$, symbolically denoted by $\Z[P\G_\infty]$,
which has an induced algebra structure. The action of $\Gn$ acts
on $P\G_n$ induces an action of the group $\G_\infty = \cup_n \Gn$
on $P\G_\infty =\cup_n P\Gn$ and thus on $\Z[P\G_\infty]$. One can
show that the orbits in $P\G_\infty$ under the action of
$\G_\infty$ give us a linear basis $\widetilde{P}_\rho$
parametrized by $\rho \in \mathcal P(\G_*)$ on $R_\Z
(P\G_\infty)$, the algebra of ${\G_\infty}$-invariants on
$\Z[P\G_\infty]$. Write $\widetilde{P}_\rho  \;
\widetilde{P}_\sigma  = \sum_{\nu} \widetilde{e}_{\rho\sigma}^\nu
\widetilde{P}_\nu.$ The homomorphisms $\theta_{n,m}$ induces a
surjective ring homomorphism $\theta_n: R_\Z (P\G_\infty)
\rightarrow R_\Z (P\Gn)$ by letting $\theta_n (\widetilde{P}_\rho)
=C_\rho(n)$ if $\Vert \rho \Vert \leq n$ and $\theta_n
(\widetilde{P}_\rho) =0$ otherwise. It follows that
$\widetilde{d}_{\rho\sigma}^\nu(n)
=\widetilde{e}_{\rho\sigma}^\nu$ for $\rho,\sigma, \nu$ such that
$\Vert\rho \Vert \leq n, \Vert \sigma\Vert \leq n$ and $\Vert
\nu\Vert \leq n$. This implies that
$\widetilde{d}_{\rho\sigma}^\nu$ is independent of $n$ and that
the stable algebra $\mathfrak A_\G$ introduced earlier is
isomorphic to the algebra $\C \otimes_\Z R_\Z (P\G_\infty)$.

A further elaboration (similar as above) using a concept of ``filling''
associated to Young diagram and the associated
convolution introduced by Ivanov-Kerov
will allow us to establish Theorem~\ref{th:interal}
for ${d}_{\rho\sigma}^\nu$.
As these are straightforward modification
of \cite{IK}, we will not repeat here and
instead refer to Sect.~8 of \cite{IK} for detail.
\subsection{Algebra generators for $\Rgn$}

For $0 \le i <n$ and $\alpha \in \Clg$, we introduce the following class function in $\Rgn$:

$$P_i(\alpha, n) = \frac{1}{ (n-i-1)!}
\cdot \mathfrak p_{-i-1}(\alpha) \mathfrak p_{-1}(1_X)^{n-i-1}\vac.$$
For $\alpha =K^c$, where $ c \in \G_*$, $P_i(\alpha, n)$ is a scalar
multiple of the characteristic function of some conjugacy class in $\Gn$.

\begin{theorem} \label{th:generator}
\begin{enumerate}
\item[{\rm (i)}] The class algebra $\Rgn$
is generated by the class functions
$\Xi^i_n(\alpha)$, where $0 \le i < n$ and $\alpha \in \G^*$.

\item[{\rm (ii)}] The class algebra $\Rgn$
is also generated by the class functions
$P^i(\alpha, n)$, where $0 \le i < n$ and $\alpha \in \G^*$.
\end{enumerate}
\end{theorem}

\begin{remark}
 Since $P_i(\alpha, n)$ depends on $\alpha$ linearly,
we may replace $\G^*$ in the theorem by any linear
basis of $\Clg$.
\end{remark}

\begin{remark}
Let us set $\G=1$ and thus $\Gn=S_n$. The power sums $\Xi_n^k$, $0\le k <n$
of the JM elements can easily be replaced
by the $k$-th elementary symmetric functions of the JM elements
in the above theorem, which is the statement of a theorem in \cite{FH}.
Note that the $k$-th elementary symmetric functions of the JM elements
is the sum of permutations in $S_n$ with exactly $n-k$ cycles.
\end{remark}

\end{document}